\documentclass[12pt]{amsart}
\usepackage{amsmath,amsfonts,amsthm,amscd,textcomp,times}
\newtheorem{theorem}{Theorem}

\newtheorem{proposition}{Proposition}
\newtheorem{lemma}{Lemma}

\newtheorem{corollary}{Corollary}
\newtheorem{remark}{Remark}
\numberwithin{equation}{section}
\numberwithin{theorem}{section}
\numberwithin{proposition}{section}
\numberwithin{lemma}{section}
\numberwithin{claim}{section}
\numberwithin{corollary}{section}
\usepackage{diagrams}
\newcommand{\bull}{\ensuremath{{}\bullet{}}}
 \newcommand{\gr}{\ensuremath{\mathbb{G}(N-n,\mathbb{C}^{N+1})}}
\newcommand{\cpn}{\ensuremath{\mathbb{P}^{N}}}

\newcommand{\dl}{\ensuremath{\partial}}

\newcommand{\ba}{\ensuremath{\begin{align*}}}
\newcommand{\ea}{\ensuremath{\end{align*}}}
\newcommand{\om}{\ensuremath{\omega}}

\newarrow{to}---->
\newarrow{mto}|--->
\newarrowtail{<=}\Rightarrow\Leftarrow\Uparrow\Downarrow
\newarrow{bboth}{<=}==={=>}
\newarrow{inject}{hooka}---{vee}
\newarrow{surject}----{>>}

\begin{document}
%%%%%%%%%%%%%%%%%%%%%%%%%%%%%%%%%%%%%%%%%%%%%%%%%%%%%%%%%%%%%%%%%%%%%%%%%%%%%%%%%%%%%%%%%%%%%%%%%%%%%%%%%%%%%%%%%%
\title[On a result of Gelfand, Kapranov and Zelevinsky]{ On a result of Gelfand, Kapranov and Zelevinsky}
\author{Sean Timothy Paul}
\email{stpaul@math.wisc.edu}
\address{Mathematics Department at the University of Wisconsin, Madison}
\thanks{Supported by an NSF DMS grant \# 0505059} 
 %%%%%%%%%%%%%%%%%%%%%%%%%%%%%%%%%%%%%%%%%%%%%%%%%%%%%%%%%%%%%%%%%%%%%%%%%%%%%%%%%%%%%%
% \title{ Projective Duality and K-Energy Asymptotics  I}
\date{August 15, 2008}
 \vspace{-5mm}
\begin{abstract}{In this paper we give new elementary proofs of basic results due to Gelfand, Kapranov, and Zelevinsky which express discriminants and results in terms of determinants of direct images of stably twisted Cayley-Koszul complexes of sheaves .}
\end{abstract}
\maketitle
\section{Statement of results}
 The main result in this paper is the following.  
 \begin{theorem} Let $X$ be a linearly normal smooth\footnote{Smoothness is only required for part $b)$. For $a)$ it is enough that $X$ be irreducible.} subvariety of $\cpn$. Let $\mathcal{V}$ be a holomorphic vector bundle over $X$. Let $({E}^{\bull}_{R}(\mathcal{V}(m)),\ \dl^{\bull}_{f})$ and $({E}^{\bull}_{\Delta}(\mathcal{V}(m)),\ \dl^{\bull}_{f})$ denote the resultant complex and the discriminant complex \footnote{Precise definitions of all the terms appearing in the statement are given in subsequent sections.} twisted by $\mathcal{V}(m)$ respectively. Then the following holds, provided $m\in \mathbb{Z}$ is sufficiently positive.
 \begin{align*}
&a) \ \mbox{{The determinant of the resultant complex is the $X$-resultant .}}\\
&\qquad \mbox{\textbf{{Tor}}}({E}^{\bull}_{R}(\mathcal{V}(m)),\ \dl^{\bull}_{f})= R^{rank(\mathcal{V})}_X(f)\ ,\quad f\in  M_{n+1, N+1}(\mathbb{C})\ .\\
\ \\
&b) \ \mbox{Assume that the dual of $X$ is non-degenerate. Then the determinant}\\
&\mbox{of the discriminant complex is the $X$-discriminant .}\\
&\qquad \mbox{\textbf{{Tor}}}({E}^{\bull}_{\Delta}(\mathcal{V}(m)),\ \dl^{\bull}_{f})= \Delta^{rank(\mathcal{V})}_X(f)\ ,\quad  f\in (\mathbb{C}^{N+1})^{\vee}\ .
\end{align*}
\end{theorem}
 
When $\mathcal{V}=\mathcal{O}_X$ this statement is due to Gelfand, Kapranov, and Zelevinsky who built upon and clarified work of Arthur Cayley.  Concerning discriminants, more general results have been obtained by Jerzy Weyman in \cite{weyman94} where he takes the higher cohomology into account. Concerning resultants, in the general case of an arbitrary \emph{coherent} sheaf $\mathcal{F}$  supported on $X$ ``Chow complexes" have been constructed by Eisenbud and Olaf-Schreyer in \cite{eisenbud}. What justifies this paper is that our proofs are completely elementary in nature. ``Elementary" means that we avoid advanced homological methods. For example, no use is made of the derived category, or the Grothendieck, Knudsen, and Mumford theory of determinants. This machinery is replaced by a very simple scaling argument, some elementary combinatorial manipulations of Chern polynomials, and a direct application of the Hirzebruch Riemann Roch formula. In this way we hope to draw the attention of geometric analysts interested in the \emph{K-stability and special metrics  problem} of K\"ahler geometry to this fascinating area of algebraic geometry and commutative algebra. 
%%%%%%%%%%%%%%%%%%%%%%%%%%%%%%%%%%%%%%%%%%

%%%%%%%%%%%%%DETERMINANTS%%%%%%%%%%%%%%%%%%%%%%%%%%%%%%%
\section{Preliminaries on Determinants}
\subsection{The Torsion of a finite dimensional chain complex}
To begin let $\left(E^{\bull},\dl^{\bull}\right)$ be a bounded complex of finite dimensional $\mathbb{C}$ vector spaces .
\[ \begin{CD} 
0@>>>E^{0}@>\dl_{0}>> \dots @>>> E^{i}@>\dl_{i}>> E^{i+1}@>\dl_{i+1}>> \dots @>\dl_{n}>> E^{n+1}@>>>0 \ . \end{CD}
\]
Recall that the \emph{determinant of the complex} $\left(E^{\bull},\dl^{\bull}\right)$ is defined to be the one dimensional vector space
\begin{align*}
  \mathbf{Det}(E^{\bull})^{(-1)^n}:=  \bigotimes_{i=0}^{n+1}(\bigwedge^{r_i}E^{i})^{(-1)^{i+1}}\quad r_i:=\mbox{dim}(E^i)\ .
\end{align*}
\begin{remark}  $\mathbf{Det}(E^{\bull})$ does not depend the boundary operators.
\end{remark}
As usual, for any vector space $V$ we set $V^{-1}:= \mbox{Hom}_{\mathbb{C}}(V,\mathbb{C})$, the dual space to
$V$.
Let $H^{i}(E^{\bull},\dl^{\bull})$ denote the $i^{th}$ cohomology group of this complex.
When $V = \mathbf{0}$ , the zero vector space , we set $\mathbf{det}(V):= \mathbb{C}$.
The determinant of the cohomology is defined in exactly the same way

\begin{align*}
\mathbf{Det}(H^{\bull}(E^{\bull},\dl^{\bull}))^{(-1)^n}:= \bigotimes_{i=0}^{n+1}(\bigwedge^{b_i}H^{i}(E^{\bull},\dl^{\bull}))^{(-1)^{i+1}}\quad b_i:= \mbox{dim}(H^{i}(E^{\bull},\dl^{\bull}))\ .
\end{align*}

\noindent We have the following well known facts (\cite{detdiv}, \cite{bgs1}). \\
 {\bf{D1}}
\emph{Assume that the complex} $\left(E^{\bull},\dl^{\bull}\right)$ \emph{is acyclic, then} $\mathbf{Det}(E^{\bull})$ \emph{is canonically trivial}
\begin{align}\label{trivial}
\tau(\dl^{\bull}):\mathbf{Det}(E^{\bull}) \cong \underline{\mathbb{C}}\ .
\end{align}
As a corollary of this we have,\\
{\bf{D2}}
 \emph{There is a canonical isomorphism \footnote{ A ``canonical isomorphism"  is one that only depends on the boundary operators, not on any choice of basis.} between the determinant of the complex and the}
\emph{determinant of its cohomology}:
\begin{align}
\tau(\dl^{\bull}):\mathbf{Det}(E^{\bull}) \cong \mathbf{Det}(H^{\bull}(E^{\bull},\dl^{\bull}))\ .
\end{align}
It is {\bf{D1}} which is relevant for our purpose. It says is that 
\emph{there is a canonically given \textbf{nonzero} element of} $\mathbf{Det}(E^{\bull})$, provided this complex is exact.
The essential ingredient in the construction of $X$-resultants (i.e Cayley-Chow forms ) and $X$-discriminants (i.e. dual varieties) consists in identifying this canonical ``section''. 
 In order to proceed we recall the \emph{Torsion} (denoted by $\mathbf{Tor}\left(E^{\bull},\dl^{\bull}\right)$) of the complex $\left(E^{\bull},\dl^{\bull}\right)$. 
 
 Define $\kappa_{i}:= \mbox{dim}(\dl_i E^{i})$, now choose $S_{i}\in \wedge^{\kappa_{i}}(E^{i})$ with $\dl_i{S_{i}}\neq 0$, then $\dl_i{S_{i}}\wedge S_{i+1}$ spans $\bigwedge^{r_{i+1}}E^{i+1}$ (since the complex is exact), that is

\begin{align*}
 \bigwedge^{r_{i+1}}E^{i+1}= \mathbb{C}\dl_i{S_{i}}\wedge S_{i+1} .
\end{align*}
With this said we define\footnote{ The purpose of the exponent $(-1)^n$ will be revealed in the next section.}
\begin{align*}
{{\mathbf{Tor}\left(E^{\bull},\dl^{\bull}\right)^{(-1)^{n}}:= (S_{0})^{-1}\otimes(\dl_0 S_{0}\wedge S_{1})\otimes (\dl_1S_{1}\wedge S_{2})^{-1}\otimes \dots \otimes (\dl_{n} S_{n})^{(-1)^{n}}}} \ .
\end{align*}
Then we have the following reformulation of {\bf{D1}}.\\ 
\noindent{\bf{D3}}
\begin{align}
 \mathbf{Tor}\left(E^{\bull},\dl^{\bull}\right) \mbox{\emph{is independent of the choices}}\ S_{i}.
 \end{align}

By fixing a basis $\{e^i_{1},e^i_{2}, \dots e^i_{r_{i}}\}$ in each of the terms $E_{i}$ we may associate to this \emph{based} exact complex a \emph{scalar}.
\begin{align*}
\mathbf{Tor}\left(E^{\bull},\dl^{\bull};\{e^i_{1},e^i_{2}, \dots, e^i_{r_{i}}\} \right)\in \mathbb{C}^*.
\end{align*}
 Which is defined through the identity:
\begin{align*}
\mathbf{Tor}\left(E^{\bull},\dl^{\bull}\right)= \mathbf{Tor}\left(E^{\bull},\dl^{\bull}; \{e^i_{1},e^i_{2}, \dots e^i_{r_{i}}\}\right)\mathbf{det}(\dots e^i_{1},e^i_{2}, \dots, e^i_{r_{i}}\dots).
\end{align*}
Where we have set
\begin{align*}
\mathbf{det}(\dots e^i_{1},e^i_{2}, \dots, e^i_{r_{i}}\dots)^{(-1)^n}:= (e^0_{1}\wedge \dots \wedge 
e^0_{r_{0}})^{-1}\otimes \dots \otimes  (e^{n+1}_{1}\wedge \dots \wedge e^{n+1}_{r_{n+1}})^{(-1)^{n}}.
\end{align*}
When we have fixed a basis of our exact complex (that is, a basis of each term in the complex) we will call $\mathbf{Tor}\left(E^{\bull},\dl^{\bull}; \{e^i_{1},e^i_{2}, \dots ,e^i_{r_{i}}\}\right)$ the Torsion of the \emph{based exact} complex. It is, as we have said, a \emph{scalar quantity}.  
\begin{remark}
In the following sections we often base the complex without mentioning it explicitly and in such cases we write (incorrectly) $\mathbf{Tor}\left(E^{\bull},\dl^{\bull}\right)$ instead of \newline
$\mathbf{Tor}\left(E^{\bull},\dl^{\bull}; \{e^i_{1},e^i_{2}, \dots ,e^i_{r_{i}}\}\right)$.
\end{remark}
 
 We have the following well known \emph{scaling behavior} of the Torsion, which we state in the next proposition. Since it is so important for us, we provide the proof, which is nothing more than the rank plus nullity theorem of linear algebra.
\begin{proposition} ( The degree of the Torsion as a polynomial in the boundary maps)
\begin{align}\label{scaling}
\begin{split}
& \mbox{deg}\mathbf{Tor}\left(E^{\bull}, \dl^{\bull}\right)= (-1)^{n+1}\sum_{i=0}^{n+1}(-1)^{i}i\mbox{\emph{dim}}(E^i) \ .
\end{split}
\end{align}
\end{proposition}
\begin{proof}
Let $\mu\in \mathbb{C}^{*}$ be a parameter. Then 
\begin{align*}
\mathbf{Tor}\left(E^{\bull},\mu \dl^{\bull}\right)^{(-1)^n}&=(S_{0})^{-1}\otimes(\mu\dl_0 S_{0}\wedge S_{1})\otimes (\mu\dl_1S_{1}\wedge S_{2})^{-1}\otimes \dots \otimes (\mu\dl_{n} S_{n})^{(-1)^{n}}\\
\ \\
&=\mu^{\kappa_0-\kappa_1+\kappa_2-\dots+(-1)^n\kappa_n}\mathbf{Tor}\left(E^{\bull}, \dl^{\bull}\right)^{(-1)^n} \ .
\end{align*}
It is clear that
\begin{align*}
\kappa_0-\kappa_1+\kappa_2-\dots+(-1)^n\kappa_n=\sum_{i=0}^{n+1}(-1)^{i+1}i(\kappa_i+\kappa_{i-1})\quad (\kappa_{n+1}=\kappa_{-1}:=0)\ .
\end{align*}
Exactness of the complex implies that we have the short exact sequence 
 \[ \begin{CD} 
0@>>>\dl_{i-1}E^{i-1}@>\iota>>E^{i}@>\dl_{i} >> \dl_iE^{i}@>>>0 \ . \end{CD}
\]
Therefore $\kappa_i+\kappa_{i-1}=\dim(E_i)$ .
\end{proof}
%%%%%%%%%%%%%%%%%%%%%%%%%%%%GEOMETRICTECHNIQUE%%%%%%%%%%%%%%%%%%%%%%%%
\section{The Geometric Technique}
 \subsection{Direct Images of Cayley-Koszul Complexes}
We let $\mathbb{C}^{\mbox{k}}$ denote the k dimensional affine space over $\mathbb{C}$.  Our concern is with irreducible subvarieties $Z$ of an affine space $\mathbb{C}^{\mbox{k}}$ associated to a smooth, linearly normal subvariety $X$ of $\cpn$. Such subvarieties $Z$ arise in the following manner.
\emph{Assume} there exists a vector subbundle $\mathcal{S}$ of the trivial bundle $\mathcal{E}:= X\times \mathbb{C}^{\mbox{k}}$ such that the image of the restriction to $I$  of the projection of $\mathcal{E}$ onto  $\mathbb{C}^{\mbox{k}}$ is $Z$, where $I$ denotes the \emph{total space} of $\mathcal{S}$. We shall always take $f$ to be a variable point in $\mathbb{C}^{\mbox{k}}$.
 
 There is the exact sequence of vector bundles on $X$
\begin{align*}
0\rightarrow \mathcal{S}\rightarrow \mathcal{E}\overset{\pi}\rightarrow \mathcal{Q}\rightarrow 0 \ .
 \end{align*}
In this case there is tautological \emph{regular} section $s$ of ${p_1}^*(\mathcal{Q})$ whose base locus is $I$.   $p_1$ denotes the projection of $\mathcal{E}$ to $X$.
We let $p_I$ denote the restriction of $p_2$ to $I$. $Z$ denotes the image of $I$ under $p_I$. This situation is pictured below in what we will call the \emph{ basic set up} following the terminology of J.Weyman (see \cite{weyman} ).
\begin{diagram}
&&{p_1}^*(\mathcal{Q})&\rTo^{\pi_2}&\mathcal{Q}\\
&& \dTo^{\pi_1}&&\dTo^{p}\\
 I&\rTo^{\iota}&{X}\times  \mathbb{C}^{\mbox{k}}&\rTo^{p_1}&X\\
\dTo ^{p_{I}}&&\dTo^{p_2}\\
  Z&\rTo ^{i}& \mathbb{C}^{\mbox{k}}&
\end{diagram}
\begin{center}\emph{In our applications we shall have that $Z$ is an irreducible algebraic {\textbf{hypersurface}} in $\mathbb{C}^{\mbox{k}}$, and that $p_I:I\rightarrow Z$ is a resolution of singularities. Therefore, in the remainder of this section we assume that $Z$ has codimension one}\end{center}

Observe that the assumption on the codimension of $Z$ in $\mathbb{C}^k$ forces $\mbox{rank}(Q)=n+1$.
In this case, following G. Kempf (see the section on ``Historical Remarks" in (\cite{kempf76}) ), we may study the irreducible equation of $Z$ (denoted by $R_{Z}(f)$ ) through an analysis of the direct image of a Cayley-Koszul complex of sheaves on ${X}\times \mathbb{C}^{\mbox{k}}$. We have the free resolution over $\mathcal{O}_{X\times \mathbb{C}^{\mbox{k}}}$
\begin{align}
(K^{\bull}(p_1^*(\mathcal{Q}^{\vee})), (s\wedge \cdot)^*)\rightarrow \iota_{*}\mathcal{O}_{I}\rightarrow 0\ ; \ K^j(p_1^*(\mathcal{Q}^{\vee})):=\bigwedge^{n+1-j}p_1^*(\mathcal{Q}^{\vee})\ .
\end{align}
More generally, let $\mathcal{V}$ denote any vector bundle on $X$. Then we will consider the \emph{twisted} 
 complex
\begin{align}
\begin{split}
&(K^{\bull}(p_1^*(\mathcal{Q}^{\vee}))\otimes p_1^*\mathcal{V} , \ (s\wedge \cdot)^*)\rightarrow \iota_{*}\mathcal{O}_{I}\otimes  p_1^*\mathcal{V}\rightarrow 0 \\
\ \\
&(s\wedge \cdot)^* \ \mbox{denotes interior multiplication} \ .
\end{split}
\end{align}
 Let $f\in  \mathbb{C}^{\mbox{k}}$, then we may pull the twisted Cayley-Koszul complex back to $X$ via the map
\begin{align}
i_{f}:X\rightarrow X\times \mathbb{C}^{\mbox{k}} \quad i_f(x):= (x,f)
\end{align}

\begin{center}\emph{Then $i_f^*(K^{\bull}(p_1^*(\mathcal{Q}^{\vee}))\otimes p_1^*\mathcal{V} , \ (s\wedge \cdot)^*)$ is an 
 {acyclic}  complex of vector bundles on $X$ whenever  $f\in \mathbb{C}^{\mbox{k}}\setminus Z$.}\end{center}

 Next we make a \emph{positivity assumption} on $\mathcal{V}$
 \begin{align}
 \mbox{We assume that}\ H^j(X, \ K^i(\mathcal{Q}^{\vee})\otimes \mathcal{V})=0 \quad \mbox{for all $i$ and all $j>0$}\ .
 \end{align}
It follows from Serre's theorem (see \cite{faisceaux}) and the Leray hypercohomology spectral sequence that the complex of \emph{finite dimensional vector spaces} is also exact 
\begin{align}
\begin{split}
&(E^{\bull}(\mathcal{V})\ , \dl_{f}):= (H^0(X, \ K^{\bull}(\mathcal{Q}^{\vee})\otimes \mathcal{V}) , \dl^{\bull}_{f})\quad f\in \mathbb{C}^{\mbox{k}}\setminus Z \\
\ \\
&\dl^{\bull}_{f}= (s(\cdot,f)\wedge \cdot)^* \ .
\end{split}
\end{align}
Choose bases $\{ e^{(\bull)}_j\}$ in each term of $E^{\bull}(\mathcal{V})$.
Then by the construction of the previous section we have a nowhere zero (and finite) \emph{rational function}\footnote{In our applications the boundary operators $\dl_f$ depend \emph{linearly} on $f$.}
\begin{align}
f\in \mathbb{C}^{\mbox{k}}\setminus Z \rightarrow   {\mathbf{Tor}}(E^{\bull}(\mathcal{V}),\dl^{\bull}_{f};\{ e^{(\bull)}_j\})\in \mathbb{C}^* \ .
\end{align}
Moreover,  an application of the Nullstellensatz to the numerator and denominator  implies at once that \emph{this rational function must be a power of} $R_{Z}(f)$. We state this in the following proposition.
\begin{proposition}
There is an integer $q$ (the {\textbf{Z-adic order}} of the determinant) such that
\begin{align}\label{ratpwr}
&{\mathbf{Tor}}(E^{\bull}(\mathcal{V}),\dl^{\bull}_{f})=R_{Z}(f)^q \ .
\end{align}
\end{proposition}
Recall that the torsion spans the determinant of the complex .
\begin{align}
\begin{split}
 &\mathbb{C}{\mathbf{Tor}}(E^{\bull}(\mathcal{V}),\dl^{\bull}_{f})=\bigotimes_{j=0}^{n+1}\bigwedge^{b_j}H^0(X,\bigwedge^j\mathcal{Q}^{\vee}\otimes \mathcal{V})^{(-1)^j} \\
\ \\
&b_j:= h^0(X, \bigwedge^{j}(\mathcal{Q}^{\vee}) \ .
\end{split}
\end{align}
 \begin{remark} In particular the determinant of the complex $(E^{\bull}(\mathcal{V}),\dl^{\bull}_{f})$ is a \textbf{polynomial}, or the reciprocal of a polynomial .
\end{remark}

If we assume that the boundary operators $\dl_f$ are \emph{\textbf{linear}} over $\mathbb{C}$ we may deduce from (\ref{ratpwr}) and (\ref{scaling}) the following well known corollary.
\begin{corollary} Assume that $\dl_f$ depends linearly on $f$. Then the degree of $Z$ can be computed as follows.
\begin{align}\label{wtdeulerchar}
q\mbox{\emph{deg}}(R_Z)= \sum_{j=0}^{n+1}(-1)^{j+1}jh^0(X, \bigwedge^{j}(\mathcal{Q}^{\vee}) \ .
\end{align}
\end{corollary}
\begin{proof}
\begin{align*}
q\mbox{\emph{deg}}(R_Z)&=(-1)^{n+1}\sum_{j=0}^{n+1}(-1)^{j}jh^0(X, \ K^{j}(\mathcal{Q}^{\vee})
\otimes \mathcal{V})\\
\ \\
&=\sum_{j=0}^{n+1}(-1)(-1)^{n+1-j}(n+1-j)h^0(X, \bigwedge^{n+1-j}(\mathcal{Q}^{\vee})
\otimes \mathcal{V})\\
\ \\
&= \sum_{j=0}^{n+1}(-1)^{j+1}jh^0(X, \bigwedge^{j}(\mathcal{Q}^{\vee})
\otimes \mathcal{V})\ .
\end{align*}
\end{proof}
\begin{corollary}\ \\
 $\mathbb{C}^k \ni f \rightarrow {\mathbf{Tor}}(E^{\bull}(\mathcal{V}),\dl^{\bull}_{f};\{ e^{(\bull)}_j\}) \in \mathbb{C}^*$
 is a constant mapping if and only if the right hand side of (\ref{wtdeulerchar}) vanishes.
\end{corollary} 
 
 %%%%%%%%%%%%%%%%%%%RESULTANTS%%%%%%%%%%%%%%%%%%%%%
 \section{Discriminants and Resultants}
\subsection{Resultants}
 Let $X \subset \cpn$ be an $n$ dimensional irreducible subvariety of $\cpn$ with degree $d$.  Then the \emph{associated hypersurface to $X$}  is defined by
\begin{align*} 
&Z_{X}:= \{L \in \mathbb{G}| L\cap X \neq \emptyset\}. \\
\ \\
&\mathbb{G}:= \mathbb{G}(N-n-1,\mathbb{P}^N) \ .
\end{align*}
It is easy to see that $Z_{X}$ is an \emph{irreducible} hypersurface (of degree $d$) in $\mathbb{G}$.
Since the homogeneous coordinate ring of the Grassmannian is a UFD, any codimension one subvariety with
degree $d$ is given by the vanishing of a section  $R_{X}$ of the homogeneous coordinate ring\footnote{See \cite{tableaux}  pg. 140 exercise 7.}
 \begin{align*}
 \{\ R_{X}=0\ \}= Z_{X}  \ ; \  R_{X} \in \mathbb{P}H^{0}(\mathbb{G},\mathcal{O}(d)).
 \end{align*}
 Following the terminology of Gelfand, Kapranov, and Zelevinsky \cite{gkz} we call $R_{X}$ the \emph{Cayley-Chow form of $X$} or simply the $X$-resultant.
Following \cite{ksz} we can be more concrete as follows.
Let $M_{n+1, N+1}^{0}(\mathbb{C})$ be the (Zariski open and dense) subspace
of the vector space of $(n+1)\times(N+1)$ matrices consisting of matrices of full rank.
We have the canonical projection
\begin{align*}
p:M_{n+1, N+1}^{0}(\mathbb{C})\rightarrow \gr ,
\end{align*}
defined  by taking the kernel of the linear transformation. This map is dominant, so the closure of the  preimage
\begin{align*}
\overline{p^{-1}(Z_{X})}\subset \overline{M_{n+1, N+1}^{0}(\mathbb{C})}=
M_{n+1, N+1}(\mathbb{C})
\end{align*}
is also an irreducible hypersurface of degree $d$ in $M_{n+1, N+1}(\mathbb{C})$.
Therefore, there is a unique\footnote{Unique up to scaling.} (symmetric multihomogeneous) polynomial (which will also be denoted by $R_{X}$)
such that 
\begin{align*}
Z:= \overline{p^{-1}(Z_{X})}= \{R_{X}(w_{ij})=0\} \ ;  \ R_{X}(w_{ij})\in \mathcal{P}^{d}[M_{n+1, N+1}(\mathbb{C})].
 \end{align*}
$R_{X}(w_{ij})$ is a polynomial of degree $d$ in the $(n+1)\times (n+1)$ \emph{minors} of $(w_{ij})$
\begin{align*}
&R_{X}(w_{ij})=\sum_{|\alpha|=d}c_{\alpha_1,\alpha_2,\dots \alpha_{b}}{P_{I_1}}^{\alpha_1}{P_{I_2}}^{\alpha_1}{\dots P_{I_{b}}}^{\alpha_b}\\
&P_I:=\mbox{det} \begin{pmatrix}w_{1i_1} & \dots & w_{1i_{n+1}}\\
       w_{2i_1} & \dots & w_{2i_{n+1}}  \\
        \dots& \dots & \dots \\
         w_{n+1i_{1}} & \dots & w_{n+1i_{n+1}}   \\ \end{pmatrix}\ .
\end{align*}
Therefore,
\begin{align}\label{scalingres}
R_{X}(\tau w_{ij})=\tau^{d(n+1)}R_{X}(w_{ij})\ .
\end{align}

Gelfand, Kapranov, and Zelevinsky (\cite{gkz}, \cite{gkzdual}, \cite{newtpoly}, \cite{gkz90}, \cite{gkz89})  have extended and clarified work of Cayley on discriminants and resultants by exhibiting $\Delta_{X}$ (respectively $R_{X}$) as the determinant of the direct image of a complex of sheaves  $(\mathcal{E}^{\bull}_{\Delta}(\mathcal{V}(m)),\ \dl^{\bull}_{f})$  (respectively $(\mathcal{E}^{\bull}_{R}(\mathcal{V}(m)),\ \dl^{\bull}_{f})$)  on $X$ through an application of the construction in the previous section. We first consider the case of resultants. In this case the basic data is chosen as follows. $\mathcal{V}$ is a holomorphic vector bundle on $X$.
\begin{align}
\begin{split}
&\mathcal{E}=X\times M_{n+1, N+1}(\mathbb{C})\\
\ \\
&\mathcal{S}=\{(x,(l_0,l_1,\dots,l_n))|\ l_i(x)=0 \ ; 0\leq i\leq n\} \\
& l_i \ \mbox{denotes a linear form on $\mathbb{C}^{N+1}$}.\\
&\mathcal{Q}\cong  \overbrace{ \mathcal{O}(1)_X\oplus \mathcal{O}(1)_X \oplus \dots \oplus \mathcal{O}(1)_X}^{n+1}  \\
\ \\
&\mathcal{V}(m)=\mathcal{V}\otimes\mathcal{O}_X(m)\quad m>>0 \in \mathbb{Z} \ .
\end{split}
\end{align}
In this situation the basic diagram takes the following shape.
\begin{diagram}
&&{p_1}^*( \mathcal{O}(1)_X \oplus \dots  \oplus \mathcal{O}(1)_X )&\rTo^{\pi_2}& \mathcal{O}(1)_X \oplus \dots  \oplus \mathcal{O}(1)_X\\
&& \dTo^{\pi_1}&&\dTo^{p}\\
 I_R&\rTo^{\iota}&{X}\times   M_{n+1, N+1}(\mathbb{C})&\rTo^{p_1}&X\\
\dTo ^{p_{I}}&&\dTo^{p_2}\\
  Z&\rTo ^{i}&  M_{n+1, N+1}(\mathbb{C})&
\end{diagram}

It is not hard to see that the direct image complex (for $m>>0$) is given as follows
\begin{align}
\begin{split}
&E_R^i(\mathcal{V}(m))=H^0(X,\ \mathcal{V}(m-(n+1-i)))\otimes \bigwedge^ {n+1-i}(\mathbb{C}^{n+1})^{\vee}\\
\ \\
&E_R^i(\mathcal{V}(m))\ni P\otimes\psi\overset{\dl_f^i}{\rightarrow} \sum_{i=0}^{n}l_iP\otimes (e_i\wedge \cdot)^*(\psi) \in E_R^{i+1}(\mathcal{V}(m))\\
\ \\
&f=(l_0,l_1,\dots,l_n)\ ,\quad \mbox{$\{\dots, e_i , \dots \}$ \emph{denotes the standard basis of} $\mathbb{C}^{n+1}$ }\\
\ \\
&\mathbb{C}\mathbf{Tor}(E_R^{\bull}(\mathcal{V}(m)), \dl^{\bull}_f)=\bigotimes^{n+1}_{j=0}\bigwedge^{P(m-j)}H^0(X,\mathcal{V}(m-j))^{(-1)^j\binom{n+1}{j}} \\
\ \\
& P(m)\ \mbox{\emph{denotes the hilbert polynomial of}} \ \mathcal{V}\ .
\end{split}
\end{align}
%%%%%%%%%%%%%%%%%%%%DISCRIMINANTS%%%%%%%%%%%%%%%%%%%
\subsection{Discriminants}
Recall that the \emph{dual variety} of $(X,L)$, usually denoted by $\widehat{X}$, is the set of hyperplanes  (identified with linear forms) $f\in \widehat{\mathbb{P}^{N}}$ (the dual projective space) tangent to $X$. Generally the dual variety is an irreducible \emph{hypersurface} in $\widehat{\cpn}$. \emph{Assuming} that $\widehat{X}$ has codimension one, let $\Delta_{X}$ denote the \emph{irreducible} polynomial defining $\widehat{X}$.   
We call $\Delta_{X}$ the \emph{discriminant} of $X$. Below we set $\widehat{d}:=\mbox{deg}(\widehat{X})$.
\begin{align*}
\widehat{X}=\{f\in  {\mathbb{P}^{N}}^{\vee}| \Delta_{X}(f)=0\}
\end{align*}
In order to study the discriminant of $X\hookrightarrow \mathbb{P}^N$ we consider the \emph{cone} over $X$, which we denote by $\tilde{X}$. Let $\{F_{\alpha}\}$ denote any generating set for the homogeneous ideal of $X$. Then
\begin{align*}
T^{1,0}(\tilde{X})=\{(p,w)\in X\times \mathbb{C}^{N+1}|\ \nabla F_{\alpha}(p)\cdot w=0 \ \mbox{for all $\alpha$} \}\overset{\iota}{\hookrightarrow} X\times \mathbb{C}^{N+1}\ .
\end{align*}
Taking the dual of the inclusion $\iota$ gives a surjection $\pi$, and hence an exact sequence
\begin{align}
0\rightarrow \mathcal{S}\rightarrow X\times (\mathbb{C}^{N+1})^{\vee}\overset{\pi}{\rightarrow} T^{1,0}(\tilde{X})^{\vee}\rightarrow 0 \ .
\end{align}
For discriminants the basic data is as follows.
\begin{align}
\begin{split}
&\mathcal{E}=X\times(\mathbb{C}^{N+1})^{\vee}\\
\ \\
&\mathcal{S}=\{(p, f) |\ T^{1,0}_p(\tilde{X})\subset \mbox{Ker}(f)\} \\
 \ \\
&\mathcal{Q}\cong T^{1,0}(\tilde{X})^{\vee}\\
\ \\
&\mathcal{V}(m)=\mathcal{V}\otimes \mathcal{O}_X(m)\quad m>>0 \in \mathbb{Z}\ .
\end{split}
\end{align}
In this situation the basic diagram takes the following shape.
\begin{diagram}
&&{p_1}^*(T^{1,0}(\tilde{X})^{\vee}  )&\rTo^{\pi_2}&T^{1,0}(\tilde{X})^{\vee}\\
&& \dTo^{\pi_1}&&\dTo^{p}\\
 I_{\Delta}&\rTo^{\iota}&{X}\times (\mathbb{C}^{N+1})^{\vee}&\rTo^{p_1}&X\\
\dTo ^{p_{I}}&&\dTo^{p_2}\\
  Z&\rTo ^{i}& (\mathbb{C}^{N+1})^{\vee} &
\end{diagram}
Observe that $T^{1,0}(\tilde{X})^{\vee}$ is isomorphic to $J_{1}(\mathcal{O}_{X}(1))$, the  sheaf of \emph{one jets} of $\mathcal{O}_{X}(1)$.   The Cayley-Koszul complex of sheaves is in this case given as follows.
\begin{align*}
&\mathcal{E}^{i}_{\Delta}(\mathcal{V}(m)):= \bigwedge^{n+1-i}J_{1}(\mathcal{O}_{X}(1))^{\vee}\otimes\mathcal{V}(m)\\
\ \\
&\dl^{i}_{f}:= (j_{1}(f)\wedge \cdot)^* \ ,\ f\in H^0(X,\mathcal{O}_{X}(1))\ .
\end{align*}
The complex of global sections is then
\begin{align}
\begin{split}
&E^i_{\Delta}(\mathcal{V}(m)):=H^0(X, \ \bigwedge^{n+1-i}J_{1}(\mathcal{O}_{X}(1))^{\vee}\otimes\mathcal{V}(m))\\
\ \\
&\dl^i_f=  (j_{1}(f)\wedge \cdot)^*\ .
\end{split}
\end{align}
Proposition (\ref{jetseq}) below together with (\ref{trivial}) implies that the torsion is given as follows.  
\begin{align}
\begin{split}
&\mathbf{Tor}(E^{\bull}_{\Delta}(\mathcal{V}(m)), \dl^{\bull}_f)=\\
\ \\
& \bigotimes^{n+1}_{j=0}\mbox{det}H^0(X, \bigwedge^{j-1}T^{1,0}_X\otimes \mathcal{V}(m-j))^{(-1)^j}
\otimes \mbox{det}H^0(X, \bigwedge^{j}T^{1,0}_X\otimes \mathcal{V}(m-j))^{(-1)^j}
\end{split}
\end{align}
 %h^0(X, \bigwedge^{i-1}T^{1,0}_X\otimes \mathcal{V}(m-i)) + h^0(X, \bigwedge^{i}T^{1,0}_X\otimes \mathcal{V}(m-i))
Then the geometric formulation of (\ref{ratpwr}) is given in the following proposition. 
\begin{proposition} (Gelfand, Kapranov, Zelevinsky \cite{gkzdual}, \cite{gkz})\\
\ \\
The complex $({E}^{\bull}_{\Delta}(\mathcal{V}(m)),\ \dl^{\bull}_{f})$ is exact provided $f$ is not tangent to $X$.\\
\ \\
The complex $({E}^{\bull}_{R}(\mathcal{V}(m)),\ \dl^{\bull}_{f})$ is exact provided $f$ does not meet $X$.
\end{proposition}
Therefore the \emph{determinants} 
\begin{align}\label{nonzero}
\begin{split}
&f\in (\mathbb{C}^{N+1})^{\vee} \rightarrow \mbox{\textbf{Tor}}({E}^{\bull}_{\Delta}(\mathcal{V}(m)),\ \dl^{\bull}_{f})\\
\ \\
&f\in M_{n+1, N+1}(\mathbb{C})  \rightarrow \mbox{\textbf{Tor}}({E}^{\bull}_{R}(\mathcal{V}(m)),\ \dl^{\bull}_{f})
\end{split}
\end{align}
are \emph{rational} functions \footnote{More, precisely, they are \emph{functions} once bases are chosen.} on $(\mathbb{C}^{N+1})^{\vee}$  (respectively $M_{n+1, N+1}(\mathbb{C})$) which are \emph{ nowhere zero} away from $\{ \Delta_{X}=0 \} $ ( respectively $\{ R_X=0\}$  ). 

\section{Proof of the Main Theorem}
The main result of this paper is the following special case of (\ref{ratpwr}). As we have mentioned already, this is due to Gelfand, Kapranov, and Zelevinsky (see \cite{gkz}) in the case when $\mathcal{V}=\mathcal{O}_X$ .\\
\ \\
\noindent \textbf{Main Theorem} (\emph{The Cayley Method})
 \begin{align*}
&a) \ \mbox{\emph{The determinant of the resultant complex is the $X$-resultant }}\\
&\qquad \mbox{\textbf{{Tor}}}({E}^{\bull}_{R}(\mathcal{V}(m)),\ \dl^{\bull}_{f})= R^{rank(\mathcal{V})}_X(f)\ ,\quad f\in  M_{n+1, N+1}(\mathbb{C})\ .\\
\ \\
&b) \ \mbox{\emph{The determinant of the discriminant complex is the $X$-discriminant}}\\
&\qquad \mbox{\textbf{{Tor}}}({E}^{\bull}_{\Delta}(\mathcal{V}(m)),\ \dl^{\bull}_{f})= \Delta^{rank(\mathcal{V})}_X(f)\ ,\quad  f\in (\mathbb{C}^{N+1})^{\vee}\ .
\end{align*}
 \emph{New Proofs}\\
\ \\
$ a)$  It follows from (\ref{nonzero}) and the Nullstellensatz that there exists $q\in \mathbb{Z}$ such that 
 \begin{align}\label{scalinga}
 \mbox{\textbf{{Tor}}}({E}^{\bull}_{R}(\mathcal{V}(m)),\ \dl^{\bull}_{f}) = R^q_X(f) \ .
  \end{align}
 Next we simply \emph{scale} the left hand side of (\ref{scalinga}), and use the fact that the boundary operators  $\dl^{\bull}_{f}$  for the resultant complex are \emph{linear} in the matrix coefficients of $f$.

On the one hand we have from (\ref{scalinga}) and (\ref{scalingres}) that
\begin{align}
 \mbox{\textbf{{Tor}}}({E}^{\bull}_{R}(\mathcal{V}(m)),\ \dl^{\bull}_{\mu f})=R^q_{X}(\mu f)=\mu^{qd(n+1)}R^q_X(f)=\mu^{qd(n+1)} \mbox{\textbf{{Tor}}}({E}^{\bull}_{R}(\mathcal{V}(m)),\ \dl^{\bull}_{ f})\ .
 \end{align}
Therefore,
\begin{align*}
\mbox{deg}(\mbox{\textbf{{Tor}}}({E}^{\bull}_{R}(m),\ \dl^{\bull}_{f}))= qd(n+1) \ .
\end{align*}

 On the other hand, linearity and (\ref{scaling}) say that we have 
 \ \\
 \begin{align*}
&  \mbox{\textbf{{Tor}}}({E}^{\bull}_{R}(\mathcal{V}(m)),\ \dl^{\bull}_{\mu f}) =  \mbox{\textbf{{Tor}}}({E}^{\bull}_{R}(\mathcal{V}(m)),\ \mu \dl^{\bull}_{f}) = \mu^{\#}\mbox{\textbf{{Tor}}}({E}^{\bull}_{R}(\mathcal{V}(m)),\ \dl^{\bull}_{f}) \\
\ \\
&\#=\sum_{i=0}^{n+1}(-1)^{i+1}i\binom{n+1}{i}\mbox{dim}(H^0(X,\ \mathcal{V}(m-i))) \ .
 \end{align*}
 It follows at once from (\ref{scaling}) that \footnote{$h^0(X,\ \mathcal{V}(m-i)))$ denotes the dimension of the vector space $H^0(X,\ \mathcal{V}(m-i)))$ .}
 \begin{align}
\begin{split}
& \mbox{deg}(\mbox{\textbf{{Tor}}}({E}^{\bull}_{R}(m),\ \dl^{\bull}_{f}))= \sum_{i=0}^{n+1}(-1)^{i+1}i\binom{n+1}{i}h^0(X,\ \mathcal{V}(m-i))
\end{split}
 \end{align}
 Next we simply evaluate the right hand side.

 In order to proceed, we recall some elementary facts about numerical functions.
 Let $f$ be any numerical function,  recall that the \emph{backwards difference} of $f$ is defined as follows.
\begin{align*}
&\Delta_- f (m):=f(m)-f(m-1)\ ; \quad  \Delta_-^{k+1}f(m):= \Delta_-^{k}f(m)-\Delta_-^{k}f(m-1) .
\end{align*} 
 Let $f_{l}(m):=m^{l}$, then it is easy to see that
\begin{align*}
 \Delta_-^{k}f_l(m)=\sum_{0\leq j \leq k}(-1)^{j}\binom{k}{j}(m-j)^{l}\ .
\end{align*}
It is not difficult to verify that
\begin{align*}
 \Delta_-^{k}f_l(m)=
\begin{cases}
k! , & \mbox{if} \ k=l \\
0,& \mbox{if} \  l<k \ .
\end{cases}
\end{align*}

We will also consider the forward difference operator.
\begin{align*}
&\Delta_+ f (m):=f(m+1)-f(m)\ ; \quad  \Delta_+^{k+1}f(m):= \Delta_+^{k}f(m+1)-\Delta_+^{k}f(m) .
\end{align*} 

\begin{align*}
 \Delta_+^{k}f_l(m)=\sum_{0\leq j \leq k}(-1)^{j+1}\binom{k}{j}(m+j)^{l}\ .
\end{align*}
Similarly we have that
\begin{align*}
 \Delta_+^{k}f_l(m)=
\begin{cases}
(-1)^{k+1}k! , & \mbox{if} \ k=l \\
0,& \mbox{if} \  l<k \ .
\end{cases}
\end{align*}
 
The dimension $h^0(X,\ \mathcal{V}(m-i)) $ is a \emph{polynomial} to which we apply the previous remarks on numerical functions.  
\begin{align*}
h^0(X,\ \mathcal{O}_X(m-i)) = \sum_{k=0}^nb_k(m-i)^k \ ;  \quad b_n=\frac{d}{n!}{rank}(\mathcal{V}) \ .
\end{align*}

\begin{align*}
 \mbox{deg}(\mbox{\textbf{{Tor}}}({E}^{\bull}_{R}(m),\ \dl^{\bull}_{f}))&= \sum_{j=0}^{n+1}(-1)^{j+1}j\binom{n+1}{j}h^0(X,\mathcal{V}(m-j)) \\
 \ \\
 &=\sum_{j=0}^{n+1}\sum_{k=0}^{n}(-1)^{j+1}j\binom{n+1}{j}b_k(m-j)^k\\
 \ \\
 &= \sum_{k=0}^{n}b_k\sum_{j=0}^{n+1}(-1)^{j}\binom{n+1}{j}(m-j)^{k+1} \\
 \ \\
 &=\sum_{k=0}^{n}b_k\Delta^{n+1}_-f_{k+1}(m) \\
 \ \\
 &=(n+1)!b_n \\
 \ \\
 &=d(n+1)rank(\mathcal{V}) \ .
 \end{align*}
Therefore,
\begin{align*}
q=rank(\mathcal{V})\ .
\end{align*}
 
Before we proceed to part $b)$ we should emphasize that in case $a)$ we needed to know \emph{apriori} that the degree of the $X$-resultant is given by
 \begin{align*}
 deg(R_X)=d(n+1) \ .
 \end{align*}
 There are many interesting formulas available which compute, in varying levels of generality, the degree of the $X$-discriminant (see \cite{kleiman}, \cite{holme}, \cite{katz} ). In the present case it is most convenient to use the following  expression for the degree of the dual variety, due to Beltrametti, Fania, and Sommese .\\
  \ \\
 \noindent {\textbf{Proposition}} (Beltrametti, Fania, Sommese \cite{bfs})\\
 \emph{The projective dual variety $\widehat{X}$ is a hypersurface if and only if $c_{n}(J_{1}(\mathcal{O}_X(1)))\neq 0$ and in this case its degree $\widehat{d}$ is given as follows}
 \begin{align}\label{degdual}
 \widehat{d}=\int_Xc_{n}(J_{1}(\mathcal{O}_X(1))) \ .
 \end{align}
 Another application of (\ref{scaling}) gives that 
 \begin{align*}
 &\mbox{deg}(\mbox{\textbf{{Tor}}}({E}^{\bull}_{\Delta}(\mathcal{V}(m)),\ \dl^{\bull}_{f}))= \sum_{j=0}^{n+1}(-1)^{j+1}jh^0(X, \bigwedge ^{j}J_{1}(\mathcal{O}_X(1))^{\vee}\otimes \mathcal{V}(m)) \ .
 \end{align*}
 The Nullstellensatz implies that there is $q\in \mathbb{Z}$ such that
 \begin{align*}
 \mbox{\textbf{{Tor}}}({E}^{\bull}_{\Delta}(\mathcal{V}(m),\ \dl^{\bull}_{f}))=\Delta_X^{q}(f) \ .
 \end{align*}
 Therefore,
 \begin{align*}
\sum_{j=0}^{n+1}(-1)^{j+1}jh^0(X, \bigwedge ^{j}J_{1}(\mathcal{O}_X(1))^{\vee}\otimes \mathcal{V}(m))  =q\widehat{d} \ .
 \end{align*}
 Our aim is to show that $q=rank(\mathcal{V})$ .

 We will require the following well known fact.
 \begin{proposition} \label{jetseq}There is an exact sequence of vector bundles on $X$.
 \begin{align}\label{xact}
 \begin{split}
 0\rightarrow \mathcal{O}_X(-1)\overset{\iota}{\rightarrow}&T^{1,0}(\tilde{X})\cong J_{1}(\mathcal{O}_X(1))^{\vee}\overset{\pi}{\rightarrow} T^{1,0}(X)\otimes \mathcal{O}_X(-1)\rightarrow 0 
 \end{split}
 \end{align}
 \end{proposition}
 Since we will need an explicit description of the maps in the sequel to this paper we recall the construction in detail.
 
The holomorphic tangent bundle to the cone on $X$ (viewed as a bundle over $X$ ) is given by
\begin{align*}
T^{1,0}(\tilde{X})=\{(p,w)\in X\times \mathbb{C}^{N+1}|\ \nabla F_{\alpha}(p)\cdot w=0 \ \mbox{for all $\alpha$} \}\overset{\iota}{\hookrightarrow} X\times \mathbb{C}^{N+1}\ .
\end{align*}
Below we abuse notation as follows. On the one hand $\pi$ denotes the map
\begin{align*}
T^{1,0}(\tilde{X})\overset{\pi}{\rightarrow} T^{1,0}(X)\otimes \mathcal{O}_X(-1)\rightarrow 0 \ .
\end{align*}
On the other hand we \emph{also} denote by $\pi$ the projection onto $\cpn$
\begin{align*}
\pi:\mathbb{C}^{N+1}\setminus \{0\}\rightarrow \cpn \ .
\end{align*}
Finally we can define $\pi$ in (\ref{xact}) by the formula (where $\pi(v)=p$ )
\begin{align*}
T^{1,0}(\tilde{X})  \ni (p,w)\rightarrow \pi(p,w):= {\pi_{*}}|_v(w)\otimes v \in T^{1,0}(X)\otimes \mathcal{O}_X(-1)\ .
\end{align*}
The rationale for this follows from the fact that for all $w\in \mathbb{C}^{N+1}$ and $\alpha \in \mathbb{C}^*$ we have 
\begin{align*}
{\pi_{*}}|_{\alpha v}(w)=\frac{1}{\alpha}{\pi_{*}}|_{v}(w) \ .
\end{align*}
To see this, one just writes down the Jacobian of 
$\pi:\mathbb{C}^{N+1}\setminus \{0\}\rightarrow \cpn $ in the affine chart $\{z_0\neq 0\}$

 \begin{align*} \pi_{*}|_z=
\begin{pmatrix}-\frac{z_1}{z^2_0} &\frac{1}{z_0}&0&0&0& \dots & 0\\ 
         -\frac{z_2}{z^2_0} &0& \frac{1}{z_0}&0&0& \dots & 0 \\
          -\frac{z_3}{z^2_0} &0&0&\frac{1}{z_0}&0& \dots & 0\\
          \dots &0&0&0 &\frac{1}{z_0}&\dots&0\\
          -\frac{z_N}{z^2_0} &0&0&0&0&\dots&\frac{1}{z_0}\\ \end{pmatrix} \ .
 \end{align*}
 $\Box$ 

  The jet exact sequence induces an exact sequence on all of the (twisted) exterior powers
 \begin{align}
  0\rightarrow \bigwedge^{i-1}T^{1,0}_X\otimes \mathcal{V}(m-i) \rightarrow \bigwedge^i J_{1}(\mathcal{O}_X(1))^{\vee} \otimes \mathcal{V}(m)\rightarrow  \bigwedge^{i}T^{1,0}_X\otimes \mathcal{V}(m-i)\rightarrow 0 \ .
  \end{align}
 Therefore, when $m>>0$ taking global sections gives the short exact sequence of vector spaces
 \begin{align}
 \begin{split}
 0\rightarrow H^0(X, \bigwedge^{i-1}T^{1,0}_X\otimes \mathcal{V}(m-i))\rightarrow &H^0(X, \bigwedge ^{i}J_{1}(\mathcal{O}_X(1))^{\vee}\otimes \mathcal{V}(m))\\
 &\rightarrow H^0(X, \bigwedge^{i}T^{1,0}_X\otimes \mathcal{V}(m-i))\rightarrow 0 \ .
 \end{split}
 \end{align}
 From which we deduce the identity
 \begin{align*} 
 h^0(X, \bigwedge ^{i}J_{1}(\mathcal{O}_X(1))^{\vee}\otimes \mathcal{V}(m))=
  h^0(X, \bigwedge^{i-1}T^{1,0}_X\otimes \mathcal{V}(m-i)) + h^0(X, \bigwedge^{i}T^{1,0}_X\otimes \mathcal{V}(m-i))\ .
   \end{align*}
Next we express the dimension of the space of global sections as an integral via the Hirzebruch-Riemann-Roch formula .
\begin{align}
 h^0(X,\bigwedge^iT^{1,0}(X)\otimes \mathcal{V}(m-i))=\int_{X}\mbox{Td}(X)\mbox{Ch}(\bigwedge^iT^{1,0}_X\otimes \mathcal{V}(m-i)) \ .
 \end{align}
 Therefore we have
 \begin{align*}
&\mbox{deg}(\mbox{\textbf{{Tor}}}({E}^{\bull}_{\Delta}(\mathcal{V}(m)),\ \dl^{\bull}_{f}))\\
\ \\
&=  \sum_{i=0}^{n+1}(-1)^{i+1}i\left(h^0(X, \bigwedge ^{i-1}T^{1,0}_X\otimes \mathcal{V}(m-i))+h^0(X, \bigwedge ^{i}T^{1,0}_X\otimes \mathcal{V}(m-i))\right)\\
\ \\
&= \int_X\mbox{Td}(X)\left(\sum_{i=0}^{n+1}(-1)^{i+1}i\{\mbox{Ch}(\bigwedge^{i-1})\mbox{Ch}(-1)+\mbox{Ch}(\bigwedge^i)\}\right)\mbox{Ch}(\mathcal{V}(m))\ .
\end{align*}
 Where we have defined $\bigwedge^i:= \bigwedge ^{i}(T^{1,0}_X(-1))$. The Chern Character expressions involving $\bigwedge^{-1}$ and $\bigwedge^{n+1}$ are taken to be zero. Next observe that
 \begin{align}
\begin{split}
& \sum_{i=0}^{n+1}(-1)^{i+1}i\{\mbox{Ch}(\bigwedge^{i-1})\mbox{Ch}(-1)+\mbox{Ch}(\bigwedge^i)\}\\
 \ \\
 &=\mbox{Ch}(-1)+\sum_{i=1}^{n}(-1)^{i+1}\mbox{Ch}(\bigwedge^i)(i-(i+1)\mbox{Ch}(-1))\\
 \ \\
 &=\sum_{i=0}^n(-1)^{i }\mbox{Ch}(\bigwedge^i)+(e^{-\omega}-1)\sum_{i=0}^n(-1)^{i+1}(i+1)\mbox{Ch}(\bigwedge^i)\ .
\end{split}
 \end{align}
 We have written $\mbox{Ch}(-1)=1+(e^{-\omega}-1)$. Next we require a combinatorial lemma .
 \begin{lemma}Let $E$ be a rank $n$ vector bundle on an $n$ dimensional variety $X$. Then the following identities hold.
 \begin{align}\label{borelserre}
\begin{split}
& \sum_{i=0}^n(-1)^i\mbox{Ch}(\bigwedge^iE)= c_{n}(E^{\vee})\\
\ \\
& \sum_{i=0}^n(-1)^{i} i\mbox{Ch}(\bigwedge^iE)=(-1)^{n}c_{n-1}(E)+ \frac{(-1)^{n+1}}{2}(c_1(E)c_{n-1}(E)-3nc_n(E))\ .
\end{split}
\end{align}
\end{lemma}
\begin{remark}
The first of these identities is due to Borel and Serre (see \cite{fulton}).
\end{remark}
\begin{proof}
 The proof of (\ref{borelserre}) is quite simple. To begin let $\{\lambda_1,\lambda_2,\dots, \lambda_n\}$ be the Chern roots of $E$. Fix $p\in \mathbb{Z}$, $0\leq p \leq n$.  We have
\begin{align*}
&p!\sum_{i=1}^{n}(-1)^ii\mbox{Ch}_p(\bigwedge^iE)=\\
\ \\
&-\sum_{1\leq i \leq n}\lambda^p_i+2\sum_{1\leq i<j\leq n}(\lambda_i+\lambda_j)^p-3\sum_{1\leq i<j<k\leq n}
 (\lambda_i+\lambda_j+ \lambda_k)^p+\dots \\
\ \\
 &+(-1)^ll\sum_{1\leq i_1<i_2<\dots <i_l\leq n}(\lambda_{i_1}+\dots+\lambda_{i_l})^p+\dots +(-1)^nn(\sum_{1\leq i \leq n}\lambda_i)^p \ .
 \end{align*}
 Fix $k\in \mathbb{N}$ and $m_j\in \mathbb{N}_+$ satisfying $\sum_{j=1}^km_j=p$. Observe that $k$ necessarily satisfies $0\leq k \leq p$.
 Then the coefficient of the monomial
 \begin{align*}
 {\lambda_{i_1}}^{m_1}{\lambda_{i_2}}^{m_2}\dots {\lambda_{i_k}}^{m_k}
 \end{align*}
   in the sum
\begin{align*}
   \sum_{1\leq i_1<i_2<\dots <i_l\leq n}(\lambda_{i_1}+\dots+\lambda_{i_l})^p
\end{align*}   
 is given by
 
 \begin{align*}
& \binom{p}{m_{1}\ m_{2}\ \dots m_{k}}\times \#\{S\ | \ S\subset \{1,2,\dots ,n\}\ ; \#(S)=l\ ; \{
{i_1},{i_2},\dots, {i_k}\}\subset S \}\\
 \ \\
 &= \binom{p}{m_{1}\ m_{2}\ \dots m_{k}}  \times\binom{n-k}{l-k} \ .
 \end{align*}
 Therefore the coefficient of ${\lambda_{i_1}}^{m_1}{\lambda_{i_2}}^{m_2}\dots {\lambda_{i_k}}^{m_k}$ 
 in the sum
\begin{align*}
p!\sum_{i=1}^{n}(-1)^ii\mbox{Ch}_p(\bigwedge^iE)
\end{align*}
 is given by
  \begin{align*}
& \binom{p}{m_{i_1}\ m_{i_2}\ \dots m_{i_k}}\sum_{k\leq l \leq n}(-1)^l l\binom{n-k}{l-k}\\
 \ \\
 &= \binom{p}{m_{i_1}\ m_{i_2}\ \dots m_{i_k}}(-1)^k\sum_{0\leq j \leq n-k}(-1)^{j}(k+j)\binom{n-k}{j} \\
 \ \\
 &= \binom{p}{m_{i_1}\ m_{i_2}\ \dots m_{i_k}}(-1)^{k+1}\Delta_{+}^{n-k}f_{1}(k)\ .
 \end{align*}
If $ p \leq n-2$ then $n-k\geq 2$. In this case the sum vanishes . Next we observe that $k$ is forced to satisfy $k=n-1$ whenever $p=n-1$. In this case $m_j=1$ for all $j$. The typical monomial is
 \begin{align*}
\frac {\lambda_{1}\lambda_2\dots\lambda_n}{\lambda_i} \ .
\end{align*}
The corresponding coefficient is then
\begin{align*}
(n-1)! (-1)^{n}\Delta_+^1f_1=(-1)^{n}(n-1)! \ .
\end{align*}
Therefore,
\begin{align*}
(n-1)!\sum_{i=1}^{n}(-1)^ii\mbox{Ch}_{n-1}(\bigwedge^iE)&=(-1)^{n}(n-1)!\sum_{i=1}^n\frac {\lambda_{1}\lambda_2\dots\lambda_n}{\lambda_i}\\
\ \\
&=(-1)^{n}(n-1)!c_{n-1}(E) \ .
\end{align*}

When $p=n=k$ there is only one monomial 
\begin{align*}
\lambda_1\lambda_2\dots \lambda_n=c_n(E) \ .
\end{align*}
The corresponding coefficient is $(-1)^{n}nn!$.

The last case is when $p=n$ and $k=n-1$.  Here the typical monomial takes the shape
\begin{align*}
\frac{\lambda_i}{\lambda_k}(\lambda_1\lambda_2\dots\lambda_n) \quad i\neq k \ .
\end{align*}
Each of which are weighted with coefficient 
\begin{align*}
\frac{1}{2}(-1)^{n+1}n! \ .
\end{align*}
The sum of such monomials is therefore
\begin{align*}
\lambda_1\lambda_2\dots \lambda_n\sum_{i=1}^n\sum_{k\neq i}\frac{\lambda_i}{\lambda_k}&=\lambda_1\lambda_2\dots \lambda_n\left(\sum_{i=1}^{n}\lambda_i\sum_{k=1}^n\frac{1}{\lambda_k}-n\right)\\
\ \\
&=\sum_{i=1}^{n}\lambda_i \sum_{k=1}^n\frac{\lambda_1\lambda_2\dots \lambda_n }{\lambda_k}-n\lambda_1\lambda_2\dots \lambda_n\\
\ \\
&=c_1(E)c_{n-1}(E)-nc_n(E)\ .
\end{align*}
Therefore,
\begin{align*}
&n!\sum_{i=1}^{n}(-1)^ii\mbox{Ch}_n(\bigwedge^iE)\\
\ \\
&=(-1)^{n}nn!c_n(E)+\frac{1}{2}(-1)^{n+1}n!\left(c_1(E)c_{n-1}(E)-nc_n(E)\right)\\
\ \\
&= \frac{1}{2}(-1)^{n+1}n! \left(c_1(E)c_{n-1}(E)-3nc_n(E)\right) \ .
\end{align*}
Justification of the first identity can be carried out in exactly the same manner,
this completes the proof of the lemma. \end{proof} 
Putting all of this together gives the identity
\begin{align*}
&\sum_{i=0}^{n+1}(-1)^{i+1}i\left(Ch(\bigwedge^{i-1})Ch(-1)+Ch(\bigwedge^{i})\right)\\
\ \\
&=c_{n}(\Omega^{1,0}_{X}(1))+\om c_{n-1}(\Omega^{1,0}_{X}(1)) \ .
\end{align*}
Therefore,
\begin{align*}
&\mbox{deg}(\mbox{\textbf{{Tor}}}({E}^{\bull}_{\Delta}(m),\ \dl^{\bull}_{f}))\\
\ \\
&=\sum_{i=0}^{n+1}\int_X Td(X)(-1)^{i+1}i\left(Ch(\bigwedge^{i-1})Ch(-1)+Ch(\bigwedge^{i})\right)Ch(\mathcal{V}(m))\\
\ \\
&= rank(\mathcal{V})\int_X(c_{n}(\Omega^{1,0}_{X})+\om c_{n-1}(\Omega^{1,0}_{X}(1)))\ .
\end{align*}
Another application of the exact sequence
\begin{align}
 \begin{split}
 0\rightarrow \mathcal{O}_X(-1)\overset{\iota}{\rightarrow}&  J_{1}(\mathcal{O}_X(1))^{\vee}\overset{\pi}{\rightarrow} T^{1,0}(X)\otimes \mathcal{O}_X(-1)\rightarrow 0 
 \end{split}
 \end{align}
shows at once that
\begin{align*}
 c_{n}(J_{1}(\mathcal{O}_X(1)))=c_{n}(\Omega^{1,0}_{X}(1))+\om c_{n-1}(\Omega^{1,0}_{X}(1)) \ .
\end{align*}

This completes the proof. $\Box$  
\begin{center}{\textbf{Acknowledgments}}\end{center}
The author is in debt to Eckart Viehweg who suggested considering the dual variety. Gang Tian, Xiuxiong Chen, Jeff Viaclovsky and Lev Borisov provided valuable discussion . This work was completed at the  Centro di Ricerca Matematica  Ennio De Giorgi in March of 2008. The author thanks Paul Gauduchon, Simon Salamon, and Gang Tian for their kind invitation. This work was supported by a NSF DMS grant 0505059 .
\bibliography{ref}
\end{document}